\documentclass[12pt]{amsart}

 \usepackage{amsmath}
\usepackage{amssymb,latexsym,cite}
\usepackage{amscd}
\usepackage{comment}
 \usepackage{xspace}
 \usepackage{verbatim}
 \usepackage{todonotes}

\newtheorem{theorem}{Theorem}[section]
\newtheorem{proposition}[theorem]{Proposition}

\newtheorem{lemma}[theorem]{Lemma}
\theoremstyle{definition}
\newtheorem{remark}[theorem]{Remark}
\newtheorem{definition}[theorem]{Definition}
\numberwithin{equation}{section}

\newcommand{\R}{{\mathbb R}}

\newcommand{\tr}{\mathrm{tr}^*}

\newcommand{\chr}[1]{\mathbf{1}\ind{#1}}

\newcommand{\BSA}{\begin{subarray}}
\newcommand{\ESA}{\end{subarray}}
\newcommand{\BAL}{\begin{aligned}}
\newcommand{\EAL}{\end{aligned}}

\newcommand{\note}[1]{\noindent\textit{#1.}\hspace{2mm}}
\newcommand{\Remark}{\note{Remark}}
\newcommand{\forevery}{\quad \forall}

\newcommand{\norm}[1]{\left \|#1\right \|}

\newcommand{\rec}[1]{\frac{1}{#1}}

\newcommand{\dist}{\mathrm{dist}\,}

\newcommand{\prt}{\partial}
\newcommand{\sms}{\setminus}

\newcommand{\ti}{\times}

\newcommand{\sbs}{\subset}

\newcommand{\ind}[1]{_{_{#1}}}

\newcommand{\sth}{such that\xspace}

\newcommand{\bvp}{boundary value problem\xspace}

\newcommand{\bdw}{\partial\Gw}

\newcommand{\qtxt}[1]{\quad\textrm{#1}}

\def\ga{\alpha}     \def\gb{\beta}       \def\gg{\gamma}
\def\gc{\chi}       \def\gd{\delta}      \def\ge{\epsilon}

      \def\gk{\kappa}      \def\gl{\lambda}
                 
    \def\gr{\rho}        
\def\gs{\sigma}       
      \def\gw{\omega}
\def\gx{\xi}                
     \def\Gd{\Delta}      

    \def\Gs{\Sigma}      
\def\Gw{\Omega}              

   \def\CM{{\mathcal M}}   
      
      \def\CC{{\mathcal C}}

\def\BBG {\mathbb G}       
   \def\BBK {\mathbb K}

\def\GTM {\mathfrak M}

\def\btr{boundary trace\xspace}
\def\tr{\mathrm{tr}}

\def\gV{{\gg V}}
\def\hr{harmonic\xspace}

\def\LVsup{$L_V$ superharmonic\xspace}
\def\LVsub{$L_V$ subharmonic\xspace}

\def\q{\quad}

\def\mbfn{\mathbf{n}}

\begin{document}

\title[Schr\"odinger equations with singular potentials]{Estimates of positive sub and super solutions of  Schr\"odinger equations with singular potentials.}
\author{Moshe Marcus}

\date{\today}

\begin{abstract}
	Consider operators  $L_{V}:=\Gd + V$ in a bounded smooth domain $\Gw\sbs \mathbb{R}^N$. Assume that $V\in C^1(\Gw)$ satisfies $V(x) \leq \bar a\,\dist(x,\bdw)^{-2}$ in $\Gw$ and that $L_V$ has a ground state $\Phi_V$ in $\Gw$. Assuming an additional condition on the behavior of $\Phi_V$ (see Section 3)
	we derive sharp, two-sided estimates of weighted integrals of positive $L_V$ harmonic functions and $L_V$ potentials. These  lead to a-priori estimates of positive $L_V$ supersolutions and subsolutions assuming (in the latter case) existence of $L_V$ \btr.
	\\ [2mm]
		[2mm]
\end{abstract}

\maketitle

\section{Introduction}
Let $\Gw$ be a bounded Lip domain in $\R^N$, $N\geq3$. We study the operator
$$L_V:=\Gd +V$$
where $V \in C^1(\Gw)$. 
We assume that the
potential $V$ satisfies the conditions:
\begin{equation}\label{Vcon1}
\exists \bar a>0\, :\quad |V(x)| \leq \bar a \gd(x)^{-2} \forevery x\in \Gw \tag{A1}
\end{equation}
$$ \gd(x)=\gd_{\bdw}(x):=\dist(x,\bdw).$$
and, 
\begin{equation}\label{Vcon2}
\gg_-<1<\gg_+. \tag{A2}
\end{equation}

Here, $\gg_+= sup\,A$ and $\gg_- =\inf\,A$ where
$$ A:= \{\gg: \int_\Gw |\nabla\phi|^2\,dx\geq \gg\int_\Gw \phi^2 V\, dx  \forevery \phi\in H^1_0(\Gw)\}.$$

Condition \eqref{Vcon1} and Hardy's inequality imply that $\gg_+>0$ and $\gg_-<0$.
If $V$ is positive then $\gg_-=-\infty$ and $\gg_+$ is \emph{the Hardy constant relative to $V$} in $\Gw$, denoted by $c_H(V)$.
If $V$ is negative, obviously $\gg_+=\infty$.

For every $\gg\in (\gg_-,\gg_+)$  there exists a Green function for  $L_{\gg V}$ in $\Gw$. 
The Green function of $L_V$ in $\Gw$ is denoted by $G_V^\Gw$.

Conditions (A1) and (A2) imply:

(i) $L_V$ has a  \textit{ground state} in the sense of Agmon \cite{Ag}. The ground state $\Phi_V$ is normalized by the condition $\Phi_V(x_0)=1$ where $x_0$ is a fixed reference point in $\Gw$.

(ii) For every $y_0\in \Gw$ and $\ge>0$ there exists a constant $C>0$ \sth 
\begin{equation}\label{Phi,G}
C^{-1}G_V(x,y_0)\leq \Phi_V(x)\leq C G_V(x,y_0) 
\forevery  x\in \Gw: \,|x-y_0|\leq \ge.
\end{equation} 	 
For a proof based on \cite{YP1}, see \cite[Lemma 1.2]{MM-Green}.

(iii) A positive \LVsup $w$ is an $L_V$ potential (i.e., it does not dominate any positive $L_V$ \hr function) if and only if $w=\BBG_V[\tau]$ for some $\tau\in \GTM_+(\Gw;\Phi_V)$ (see \cite{An-SLN}). Here  $\GTM(\Gw;\Phi_V)$ denotes the space of Borel measures $\tau$ \sth $\int_\Gw \Phi_V d|\tau|<\infty$.

$L_V$ is weakly coercive in the sense of Ancona. (A proof, due to \cite{Pincho-comm}, is provided  in \cite[Lemma 1.1]{MM-Green}.)
Therefore, by Ancona \cite{An87}: \\ [1mm]
$\bullet$ $L_V$  possesses a Martin kernel $K_V$ such that, for every $y\in \bdw$, $x\mapsto K_V(x,y)$ is positive $L_V$ harmonic in $\Gw$ and vanishes on $\bdw\sms \{y\}$ and the Representation Theorem holds:

\textit{If $u$ is a positive $L_V$-\hr function then there exists $\nu\in \GTM_+(\bdw)$ ($=$ the space of positive, bounded Borel measures) \sth}
	\begin{equation}\label{rep}
	u(x)=\int_{\bdw} K_V(x,y) d\nu(y)=:\BBK_V[\nu] \q x\in \Gw.
	\end{equation}

$\bullet$ The \textsc{Boundary Harnack Principle} (briefly BHP) holds. (See also Ancona \cite{An_MZ}.)\\

The present paper is devoted to the derivation
of weighted integral estimates of positive \LVsup and \LVsub functions.
The weight $W$ is given by,
\begin{equation}\label{weight}
W:=\frac{\Phi_V}{\Phi_0}.
\end{equation}
The  estimates are sharp and two sided (see Theorem \ref{LVsup} below). In this sense the weight $W$ is optimal. The derivation is based on assumptions (A1), (A2) and an additional condition on the behavior of the ground state, (see (C1) in section 3). 

Such estimates have been derived in \cite{MT2} and \cite{MT3} when $V$ is the Hardy potential or $V$ is of the form
\begin{equation}\label{S1.1}\BAL
V=\gg V_F,\q 
V_F= \frac{1}{\gd_F^2},\q  \gd_F(x)=\dist(x,F)
\EAL\end{equation} 
where  $F\sbs \bdw$ is a smooth k-dimensional manifold without boundary. In these cases a sharp two sided estimate of $\Phi_V$ is available. This fact was crucial for  ) The estimates have been applied in the study of positive solutions of a family of semilinear boundary value problems involving $L_V$ and a nonlinear term.

\vskip 3mm

Linear and nonlinear boundary value problems for operators $L_V$, with $V$ as in \eqref{S1.1},  have been investigated by many authors. The cases $F=\bdw$ or $F$ a singleton have been most frequently investigated. Following is a list of some recent works in the area:\\ [1mm] 

\noindent Bandle, Moroz and Reichel \cite{BMR1}, \cite{BMR2}, Marcus and P.T. Nguyen \cite{MT2}, \cite{MT3},
Gkikas and Veron \cite{GkVe}, P.T. Nguyen \cite{Ng}, Y. Du and L. Wei \cite{DW2015}, \cite{DW2017},  Marcus and Moroz \cite{MM-VM}, Chen and Veron \cite{HC-LV}, Gkikas and Nguyen \cite{Gk_Ng1}, \cite{Gk_Ng2}.

The main tools used in the paper are potential theoretic results (mentioned above) and estimates of the
Green and Martin kernels \cite{MM-Green}. 

The plan of the paper: Section 2 is devoted to notations and statement of some results from the literature. The main results - two sided estimates of $\BBK_V[\nu]$, $\nu\in \GTM_+(\bdw)$ and $\BBG_V[\tau]$, $\tau \in \GTM_+(\Gw;\Phi_V)$ - are stated in Section 3 and proved in Sections 4 and 5. In order to simplify the presentation we assume that $\Gw$ is a smooth domain.\footnote{Employing local coordinates as described in Section 2, the main results can be extended to Lipschitz domains.} In section 6 we describe a family of potentials $V$ that satisfy conditions (A1), (A2), (C1) and include the potentials studied in \cite{MT2} and \cite{MT3}. Finally in Section 7 we discuss the notion of $L_V$ \btr and apply the previous estimates to derive estimates of positive $L_V$ superharmonic and $L_V$ subharmonic functions. In the first case we use the Riesz representation formula. In the second case we show that, assuming existence of the $L_V$ trace, a similar representation formula holds.

\section{Notation and preliminaries}
Denote,
$$T(r,\gr)= \{\gx= (\gx_1, \gx')\in \R\ti\R^{N-1}:    |\gx_1|<\gr,\;  |\gx'|<r\}.$$
Assuming that $\Gw$ is a bounded Lipschitz domain, there exist positive numbers  $r_0$ , $\gk$ \sth, for every $y\in \bdw$, there exist: (i) a  set of Euclidean coordinates $\gx=\gx_y$ centered at $y$ with the positive $\gx_1$ axis pointing in the direction of $\mbfn_y$\footnote[1]{If $\Gw$ is smooth, $\mbfn_y$ denotes the inward normal at $y$. If $\Gw$ is Lipschitz, $\mbfn_y$ denotes an approximate normal.} and (ii) a function $F_y$ uniformly Lipschitz in $\R^{N-1}$ with Lipschitz constant $\leq \kappa$ \sth
\begin{equation}\label{ngh_y}\BAL
Q_y(r_0,\gr_0):=\, &\Gw\cap T_y(r_0,\gr_0)\\
=\, &\{ \gx= (\gx_1, \gx'): F_y(\gx')<\gx_1<\gr_0,\; |\gx'|<r_0\},
\EAL\end{equation}
where $ T_y(r_0,\gr_0)=y  +  T(r_0,\gr_0)$ in coordinates $\gx=\gx_y$ and $\gr_0 = 10\gk r_0.$
Without loss of generality, we assume that $\kappa>1$.
The set of coordinates $\gx_y$ is called a standard set of coordinates at $y$ and $T_y(r,\gr)$ with $0<r\leq r_0$ and $\kappa r< \gr\leq 10\kappa r$ is called a standard cylinder at $y$.

If $\Gw$ is a bounded $C^2$ domain  there exists $\bar \gb>0$ \sth  for every $x\in \Gw_{\bar \gb}$ there is a unique point $\gs(x)\in \bdw$ \sth
$$|x-\gs(x)|=\gd(x)$$
and $x\mapsto \gd(x)$ is in $C^2(\Gw_{\bar \gb})$ while $x\mapsto \gs(x)$ is in $C^1(\Gw_{\bar \gb})$.
 The set of coordinates $(\gd,\gs)$ defined in this way in $\Gw_{\bar \gb}$
 is called \textit{the flow coordinates set}.
  We denote
\begin{equation}\label{notation1}\BAL
D_\gb=\{x\in \Gw:\gd(x)>\gb\}, &\quad \Gw_\gb=\{x\in \Gw:\gd(x)<\gb\},\\
\Gs_\gb=\{x\in \Gw:&\gd(x)=\gb\}.
\EAL\end{equation}
In the sequel we denote
$$\gb_0 :=\min (r_0, \bar \gb).$$

\noindent\textit{Notation.}\hskip 2mm  Let $f_i$, $i=1,2$, be positive functions on some domain $X$. Then the notation
$f_1\sim f_2$ in $X$
means: there exists  $C>0$ \sth
$$\rec{C}f_1\leq f_2\leq C f_1 \qtxt{in }X.$$
The notation $f_1\lesssim f_2$ means: there exists  $C>0$ \sth
$f_1\leq Cf_2$ in $X$. The constant $C$ will be called a \textit{similarity constant}.
\vskip 2mm

We state below the Boundary Harnack Principle in  Lipschitz domains, due to Ancona [An87]:

      \begin{theorem} \label{BHP}
      	
      	Let $y\in \bdw$ and let $\gx_y$ and $T_y(r,\rho)$ be a standard set of coordinates and a standard cylinder at $y$. 
      Let $\gw:=\Gw\cap T_y(r,\rho)$ and $A=(0,...,0,\gr/2)$, $A'=(0,...,0,3\gr/4)$ in the coordinates $\gx=\gx_y$.
      
	There exists a constant $c$ depending only on $N,\, \bar a$ and ${ \frac {\rho } { r}}$ such that whenever $u$ is a positive $L_V$ harmonic function  in $\gw$  that vanishes continuously in $\prt \gw  \cap T_y(r, \rho )$
	the following inequality holds:
	\begin{equation}\label{u/v+}
	 c^{-1}r^{N-2}\, G_{V}(x,A')\leq \frac{u(x)}{u(A)}\leq c\, r^{N-2} \, G _{V}(x,A'),  
	\end{equation}
	for every $x\in \gw  \cap T_y({ \frac {  r} {2 }}; { \frac {   \rho } {2 }})$.
	In particular,  for any pair $u,v $ of  positive $L_V$ harmonic functions in  $ \gw $ that vanish on $\partial \gw  \cap T_y(r, \rho )$:
	\begin{equation}\label{u/v}
	u(x)/v(x)\le Cu(A)/v(A),  \forevery x\in \gw \cap T_y(r/2, \rho /2))
	\end{equation}
	where $C=c^2$.
\end{theorem}

The following is a well-known consequence of BHP (see e.g.   \cite[Lemma 3.5]{An_MZ}).
Here $x_0$ is a reference point in $\Omega \setminus \overline  T(r, \rho )$ and $\mbfn_y$ is the unit vector in the direction of the $(\gx_y)_1$ axis.

\begin{proposition}\label{c:BHP}
	There exists a constant $C$  such that for all  $x = y+t\mbfn_y$,  $|t|  \leq { \frac {  3} {4 }} \rho $, \begin{equation}
	C^{-1}t^{2-N}\leq  K_V(x, y) G_V(x,x_0)\leq C\,t^{2-N}.
	\end{equation}
	and $C$ can be chosen depending only on $\bar a$, ${ \frac {   \rho } {r }}$ and $N$.
\end{proposition}

Estimates of the Green and Martin kernels \cite{MM-Green} will be frequently used in the sequel. These are valid in \emph{bounded Lipschitz domains}. 

\vskip 3mm

\begin{theorem}\label{G-estA} Assume (A1), (A2) and $N\geq 3$.
	
	Then, for every $b>0$ there exists
	a constant $C(b)$, depending also on $N, r_0, \kappa, \bar a$, \sth: if $x,y\in \Gw$ and
	\begin{equation}\label{d>>r}
	|x-y|\leq \rec{b}\min(\gd(x),\gd(y))
	\end{equation}
	then
	\begin{equation}\label{G-est.1}\BAL
	\rec{C(b)}|x-y|^{2-N}\leq G_V(x,y)\leq C(b)|x-y|^{2-N}.
	\EAL\end{equation}
\end{theorem}
\medskip

In the next theorems, $C$ stands for a constant depending only on $r_0,\kappa, \bar a$ and $N$.

\begin{theorem}\label{G-estB} Assume (A1), (A2) and $N\geq 3$.
	
	If $x,y\in \Gw$ and
	\begin{equation}\label{d<r0}
	\max(\gd(x),\gd(y))\leq r_0/10\kappa
	\end{equation}
	\begin{equation}\label{d<<r}
	\min(\gd(x),\gd(y))\leq \frac{|x-y|}{16(1+\kappa)^2}
	\end{equation}
	then
	\begin{equation}\label{G-est.2}\BAL
	\rec{C}|x-y|^{2-N}\frac{\Phi_V(x)\Phi_V(y)}{\Phi_V( x_y)^2}&\leq G_V(x,y)\\
	&\leq
	C|x-y|^{2-N}\frac{\Phi_V(x)\Phi_V(y)}{\Phi_V( x_y)^2}.
	\EAL\end{equation}
	The point $x_y$ depends on the \emph{pair} $(x,y)$. If
	$$\hat r(x,y):=|x-y|\vee \gd(x)\vee \gd(y)\leq r_0/10\kappa$$
	$x_y$ can be chosen arbitrarily in the set
	\begin{equation}\label{Axy}
	A(x,y):=\{z\in \Gw: \rec{2}\hat r(x,y)\leq \gd(z)\leq 2\hat r(x,y)\} \cap B_{4\hat r(x,y)}(\frac{x+y}{2})\}.
	\end{equation}
	Otherwise set $ x_y= x_0$
	where $x_0$ is a fixed reference point.
\end{theorem}

\begin{theorem}\label{K-est} Assume (A1), (A2) and $N\geq 3$.
	
	If $x\in \Gw$, $y\in \bdw$ and $|x-y|< \frac{r_0}{10\kappa}$ then
	\begin{equation}\label{K-est.1}\BAL
	\rec{C} \frac{\Phi_V(x)}{\Phi_V(x_y)^2}|x-y|^{2-N}\leq K_\Gw^\gV(x,y)\leq
	C \frac{\Phi_V(x)}{\Phi_V(x_y)^2}|x-y|^{2-N},
	\EAL\end{equation}
	where $x_y$ is an arbitrary point in $A(x,y)$.
\end{theorem}
	

\section{Main results}

In the results stated below $\Gw$ is a bounded $C^2$ domain in $\R^N$. The results will be extended,  in a separate note, to the case of bounded Lipschits domains.

The first result provides a sharp estimate of positive $L_V$ harmonic functions.
\begin{theorem}\label{Knu}     Assume conditions (A1), (A2). 
	In addition assume that, for every $a>a_0>>1$ and every pair $x,z\in \Gw_{\gb_0}$ \sth $z$ lies on the normal to $\bdw$ at $\gs(x)$ (= nearest point to $x$ on $\bdw$): 
	\[
a\gd(x)\leq \gd(z)\Longrightarrow
	\frac{\Psi_V(x)}{\Psi_V(z)}\leq c(a) \frac{\gd(x)^{\ga^*}}{\gd(z)^\ga},  \tag{C1}
	\]
where
\begin{equation}\label{ga*}
   0\leq \ga-\ga^*<1/2.
\end{equation}	
	 
	Then
	\begin{equation}\label{KV1}
	\rec{C}\|\nu\|\leq\int_{\Gw}\frac{\Phi_V}{\Phi_0}\BBK_V[\nu]dx\leq C\|\nu\| \forevery \nu\in \CM_+(\bdw),
	\end{equation}
	where $C$ depends on $\bar a,\Gw$ and the constants in (C1).
\end{theorem}	

\remark\label{rem01}	(i) The lower estimate requires only conditions (A1) and (A2). Condition (C1) is used in the derivation of the upper estimate when $a\gd(x)\leq |x-y|$.

(ii) If \eqref{ga*} holds for $\ga$ and $\ga^*$  then, it also holds with  $\ga+\ge$, $\ga^*$ with $\ge>0$ \sth 
$\ge+\ga-\ga^*< 1/2$. Therefore, without loss of generality, we may assume that $\ga\neq 1/2$.\\

In the following two theorems we present  estimates of $L_V$ potentials. Recall that $w$ is an $L_V$ potential iff $w=\BBG_V[\tau]$ for a positive measure $\tau\in \CM(\Gw;\Phi_V)$.

\begin{theorem}\label{G-I} Assume (A1) and (A2).
	Then there exists a  constant $c$ depending on $\bar a$, $r_0$ and $\kappa$ \sth, for every $\tau\in \GTM_+(\Gw;\Phi_V)$,
	\begin{equation}\label{G-I1}
	\rec{c}\int_{\Gw_{\gb_0/4}} \Phi_V d\tau\leq \int_{\Gw_{\gb_1}}\frac{\Phi_V}{\Phi_0}\BBG_V[\tau]dx
	\end{equation}
	where $\gb_1= \gb_0/48(2+\kappa)^2$ and
		\begin{equation}\label{G-I1'}
	\rec{c}\int_{\Gw} \Phi_V d\tau\leq \int_{\Gw}\frac{\Phi_V}{\Phi_0}\BBG_V[\tau]dx
	\end{equation}
	 
\end{theorem}

\begin{theorem}\label{G-II}
	 Assume (A1), (A2) and (C1).

Then there exists $c'>0$, depending on $\bar a$, $\Gw$ and the constants in (C1) \sth  for every $\tau\in \GTM_+(\Gw;\Phi_V)$
	\begin{equation}\label{G-I2}
\int_{\Gw}\frac{\Phi_V}{\Phi_0}\BBG_V[\tau] dx \leq c'\int_\Gw \Phi_V d\tau.
\end{equation}
\end{theorem}

\remark\label{rem02}\hskip 3mm (i) As before, without loss of generality we may assume that in (C1), $\ga\neq 1/2$ (see Remark \ref{rem01}).

(ii)  See also Lemmas \ref{L-I1} and  \ref{L-I2} below for estimates of surface integrals on $\Gs_\gb= \{x\in\Gw: \gd(x)=\gb\}$,  $\gb<\gb_0$.\\

Estimates as above and a version of Theorem \ref{LVsup} have been proved in \cite{MT2} for $V=\gg/\gd^2$ and in \cite{MT3} for $V=\gg V_k$ where $V_k=\gd_{F_k}^{-2}$ and $F_k$ is a smooth $k$-dimensional manifold without boundary. The estimates
in \cite{MT3} required $\gg < \min(c_H(V_k),\rec{4}(2(N-k)-1))$. The present estimates apply to a family of potentials which include those mentioned above (see Section 6) and  require only $\gg<c_H(V_k)$.

In \cite{MT2}, \cite{MT3} and \cite{MM-VM} the estimates have been applied to a study of semilinear \bvp with absorption nonlinearity. (In those papers the definition of \btr was different from the $L_V$ trace used in Section \ref{HM}. However, the two definitions are equivalent in the specific cases studied there.)  

\vskip 3mm

\section{Estimates of $L_V$ harmonic functions}

\noindent\textbf{Proof of Theorem \ref{Knu}.} 
Let $y\in\bdw$ and $b>1$. Put $$\CC_b(y)=\{x\in \Gw:|x-y|\leq b\gd(x),\; \gd(x)<\ge_b\}$$
where $0<\ge_b<\gb_0$ and $b$  are chosen so that
$$\CC_{2b}(y)\sbs \Gw_{\gb_0}\cup\{y\} \forevery y\in \bdw.$$

By Proposition \ref{c:BHP} there exists $t_0\in (0, \gb_0)$ \sth, for every $y\in \bdw$,
\begin{equation}\label{mk1}
K_V(y+t\mbfn_y, y) G_V(y+t\mbfn_y, x_0)\sim \, t^{2-N}  \q t\in (0,t_0)
\end{equation}
with similarity constant dependent on $\bar a, \Gw$ but independent of $y$.

Let $y\in \bdw$ and $b>1$. Taking $t_0$ sufficiently small (depending on $b$)
$$\CC_b(y):=\{\langle x-y, \mbfn_y\rangle >0, \;|x-y|\leq b\gd(x),\; \gd(x)<t_0\} \sbs \Gw$$
for every $y\in \bdw$.

 By \eqref{mk1} and the strong Harnack inequality,
$$\frac{K_0(x,y)}{K_V(x,y)} \sim \frac{G_V(x,x_0)}{G_0(x,x_0)}\sim\frac{\Phi_V(x)}{\Phi_0(x)} \forevery x\in \CC_b(y)\cap \Gs_\gb, \forevery \gb \in (0,t_0).$$ 
Hence,
\begin{equation}\label{mk0}
\int_{\CC_b(y)\cap \Gs_\gb}\frac{\Phi_V(x)}{\Phi_0(x)} K_V(x,y)dS_x \sim \int_{\CC_b(y)\cap \Gs_\gb} K_0(x,y)dS_x,
\end{equation}
with similarity constants independent of $y\in \bdw$ and $ \gb\in(0,t_0)$.
Hence there exist positive constants $c_1$ and $c_2(b)$ independent of $y\in \bdw$ and $\gb\in (0,t_0)$ \sth,

\begin{equation}\label{mK2}
\rec{c_1}\leq \int_{ \Gs_\gb}\frac{\Phi_V(x)}{\Phi_0(x)} K_V(x,y)dS_x 
\end{equation}
and 
\begin{equation}\label{mK3}
\int_{\CC_b(y)\cap \Gs_\gb}\frac{\Phi_V(x)}{\Phi_0(x)} K_V(x,y)dS_x \leq c_2(b).
\end{equation}

 By Fubini's theorem, the lower estimate in \eqref{ga*} is a consequence of \eqref{mK2}.

Next we show that there exists $\gl>-1$ \sth
\begin{equation}\label{mK3.1}
\int_{\Gs_\gb\sms\CC_b(y)}\frac{\Phi_V(x)}{\Phi_0(x)} K_V(x,y)dS_x \leq c_3(b) \gb^\gl   \forevery \gb\in (0,t_0)
\end{equation}
where $c_3(b)>0$ is independent of $\gb$ and $y$. 

By Theorem \ref{K-est},
\begin{equation}\label{e:K-est}
\frac{\Phi_V(x)}{\Phi_0(x)}K_V(x,y)\leq C \frac{\Phi_V^2(x)}{\Phi_V^2(x_y)}|x-y|^{2-N}\gd(x)^{-1}
\end{equation}
where, if $|x-y|\leq \gb_0$, we choose $x_y=(\gd(x)+|x-y|)\mbfn_{\gs(x)}$. Thus $x$ and $x_y$ are on the normal to $\bdw$ at $\gs(x)$. (Recall that $\gs(x)$ denotes the closest point to $x$ on $\bdw$.)
Consequently, by \eqref{e:K-est} and (C1) (with $z=x_y$) if 
$$x\in E_y(\gb,t_0):= 
\{\gd(x)=\gb: b\gb<|x-y|\leq t_0\},$$
then
$$\frac{\Phi_V(x)}{\Phi_0(x)}K_V(x,y)\leq C \gb^{2\ga^*-1}
|x-y|^{2-N-2\ga}.$$
Hence, 
$$\int_{E_y(\gb,r_0)}\frac{\Phi_V(x)}{\Phi_0(x)} K_V(x,y)dS_x \leq C \gb^{2\ga^*-1}\int_{b\gb}^{\gb_0} s^{-2\ga}ds \leq C' \gb^{2(\ga^*-\ga)}. $$
Here we assumed that $\ga\neq 1/2$. As mentioned in Remark \ref{rem01}, this does not involve a loss of generality in Theorem \ref{Knu}.

On the other hand, if $|x-y|\geq t_0$, \eqref{e:K-est}
holds with $x_y=x_0$ so that $\Phi_V(x_y)=1$. Therefore, by \eqref{e:K-est},

$$ \int_{\Gs_\gb\cap[|x-y|\geq r_0] }\frac{\Phi_V(x)}{\Phi_0(x)} K_V(x,y)dS_x \leq C.$$
The last two inequalities imply \eqref{mK3.1} with $\gl:=2(\ga^*-\ga)$. By \eqref{ga*}, $-1<\gl\leq 0$.

By \eqref{mK3} and \eqref{mK3.1},
\begin{equation}\label{mK_Gs_gb}
\int_{\Gs_\gb} \frac{\Phi_V(x)}{\Phi_0(x)} K_V(x,y)dS_x \leq C(b)\gb^\gl.
\end{equation}
Combining this inequality with \eqref{mK2} and
integrating over $\gb\in (0,\gb_0)$ we obtain 
$$ \rec{C}\leq\int_{\Gw_{\gb_0}}\frac{\Phi_V(x)}{\Phi_0(x)} K_V(x,y)dS_x \leq C.$$
Assuming that $\gd(x_0)>\gb_0$, the integral over $\Gw\sms \Gw_{\gb_0}$ is also bounded by a constant independent of $y$. (Recall that $K_V(x_0,y)=1$ for every $y\in \bdw$ and use Harnack's inequality for positive $L_V$ harmonic functions.) Therefore this inequality implies \eqref{KV1} in the case $\nu=\gd_y$. By Fubini's theorem, \eqref{KV1} follows in the general case.
\qed

\remark\label{rem03} Inequality \eqref{mK_Gs_gb} is of interest in itself. It should be noted that this inequality is valid even if $\ga=1/2$ but in that case $\gl=2(\ga^*-\ga) -\ge$ where $\ge>0$ satisfies $\ga+\ge-\ge^* < 1/2$.

\section{Estimates of $L_V$ potentials}

In this section we prove Theorems \ref{G-I} and \ref{G-II}. The proofs are based on two lemmas.

\begin{lemma}\label{L-I1} Assume that (A1) and (A2) hold. Let $\tau\in \GTM_+(\Gw;\Phi_V)$ and denote
	$$I_1(\gb):=\rec{\gb}\int_{\Gs_\gb}\Phi_V(x)\int_\Gw G_V(x,y) \gc_{a\gb}(|x-y|)d\tau(y) dS_x$$
	where $\gc_{s}(t)=\chr{(0,s)}(t)$ and $a\geq 16(\kappa+2)^2$.	Then there exists a constant $c$ depending only on $a$, $\bar a$ and $\Gw$ \sth,
	\begin{equation}\label{I1}
	\rec{c}\int_{\Gw_{3a\gb/2}}\Phi_V d\tau \leq I_1(\gb)\leq c\int_\Gw \Phi_V d\tau \forevery \gb\in (0, \gb_0/3a).
	\end{equation}
\end{lemma}

\proof
 The domain of integration in $I_1(\gb)$ is $\{(x,y)\in \Gs_\gb\ti \Gw: |x-y|<a\gb\}$. We partition the domain of integration  into three parts and estimate each of the resulting integrals separately. Accordingly we denote:
\begin{equation*}\BAL
I_{1,1}(\gb)&:= \rec{\gb}\int_{\Gs_\gb}\Phi_V(x)\int_{\gb/a\le \gd(y)\le\gb} G_V(x,y) \gc_{a\gb}(|x-y|)d\tau(y) dS_x,\\
I_{1,2}(\gb) &:= \rec{\gb}\int_{\Gs_\gb}\Phi_V(x)\int_{\gd(y)\le \gb/a} G_V(x,y)
\gc_{a\gb}(|x-y|)d\tau(y) dS_x,\\
I_{1,3}(\gb)&:= \rec{\gb}\int_{\Gs_\gb}\Phi_V(x)\int_{\gb\le \gd(y)} G_V(x,y) \gc_{a\gb}(|x-y|)d\tau(y) dS_x
\EAL\end{equation*}
so that $I_1=I_{1,1} +I_{1,2} +I_{1,3}$.
\vskip 2mm

\noindent\textit{Estimate of $I_{1,1}(\gb)$.}\\ [2mm]
By  the  Hardy (chain) inequality (see e.g. \cite[Lemma 3.2]{MM-Green}) , there exists $C(a)>0$ \sth, if
\[\gb/a\le \gd(y)\leq \gb, \q x\in \Gs_\gb,\q  |x-y|\leq a\gb  \tag{*}\]
then
\begin{equation}\label{Har1}
\rec{C(a)}\Phi_V(x)\leq \Phi_V(y)\leq C(a)\Phi_V(x).
\end{equation}

By Theorem \ref{G-estA},  if (*) holds
then
$$\rec{c}|x-y|^{2-N}\le G_V(x,y)\le c |x-y|^{2-N}$$
for some constant $c=c(a)$.
Hence,
\begin{equation}\label{temp2}\BAL
I_{1,1}(\gb)&\sim \rec{\gb} \int_{\Gs_\gb}\int_{\gb/a\le \gd(y)\le\gb} |x-y|^{2-N}\gc_{a\gb}(|x-y|) \Phi_V(y) d\tau(y)dS_x \\
&=\rec{\gb}\int_{\gb/a\le \gd(y)\le\gb} \int_{\Gs_\gb}|x-y|^{2-N}\gc_{a\gb}(|x-y|)dS_x \Phi_V(y) d\tau(y)
\EAL\end{equation}
For every $y$ \sth $\gb/a\le \gd(y)\le\gb$, 
the domain of integration contains the set $\{x\in\Gs_\gb: \gb -\gb/a\leq |x-y|\leq a\gb\}$.  Therefore, for $y$ as above,
\begin{equation}\label{temp2'}
\BAL 
(a-1)\gb \lesssim\int_{\Gs_\gb}|x-y|^{2-N}\gc_{a\gb}(|x-y|)dS_x\lesssim\int_0^{a\gb} dr=a\gb,
\EAL\end{equation}
Hence by \eqref{temp2},
\begin{equation}\label{est2.3}
I_{1,1}(\gb)\sim  \int_{{\gb/a\le \gd(y)\le\gb}}\Phi_V d\tau
\end{equation}
with similarity constant depending on $a$ and $\bdw$.\\ [2mm]

\noindent\textit{Estimate of $I_{1,2}(\gb)$.}\hskip 2mm Here we assume that $\gb<\gb_0/3a$.
Since $\gd(y)< \gb/a $ it follows that in the domain of integration of $I_{1,2}$,
\begin{equation}\label{interval1}
a\gb\geq |x-y|\geq \gb-\gb/a>\gd(y)(a-1)
\end{equation}
Thus the pair $(x,y)$ satisfies the conditions of
Theorem \ref{G-estB} and consequently,
\begin{equation}\label{est.G3}
\Phi_V(x)G_V(x,y) \sim \frac{\Phi_V(x)^2}{\Phi_V(x_y)^2}\Phi_V(y)|x-y|^{2-N},
\end{equation}
where $x_y$ may be chosen as follows:
$x_y:=\eta + |x-y|\mathbf{n}_\eta$ with $\eta\in \bdw$  the  closest point to $x$. 

Then $\gd(x_y)=|x-y|$ and, by \eqref{interval1}, 
$$\gb(1-\rec{a})\leq \gd(x_y)\leq a\gb, \q |x-x_y|\leq |x-y|\leq a\gb.$$
 Hence, by the Hardy (chain) inequality (see \cite[Lemma 3.2]{MM-Green}), there exists $c'(a)>0$ \sth
$$\rec{c'}\Phi_V(x)\leq \Phi_V(x_y)\leq c'\,\Phi_V(x).$$
Therefore, by \eqref{est.G3},
\begin{equation}\label{temp1}
\BAL \Phi_V(x)G_V(x,y)\sim  |x-y|^{2-N}\Phi_V(y).
\EAL
\end{equation}
with similarity constant depending on $a$.
Consequently  (using Fubini's theorem) we obtain 
\begin{equation}\label{temp3}\BAL
I_{1,2}(\gb) \sim \rec{\gb}\int_{\gd(y)\leq\gb/a} \int_{\Gs_\gb\cap[|x-y|\leq a\gb]}|x-y|^{2-N} dS_x \Phi_V(y) d\tau(y).
\EAL\end{equation}
Hence,
\begin{equation}\label{est2.3'}
I_{1,2}(\gb)\sim  \int_{\gd(y)\leq \gb/a}  \Phi_V(y) d\tau(y).
\end{equation}

\noindent\textit{Estimate of $I_{1,3}(\gb)$.}
 In this case, as $|x-y|<a\gb$, $\gd(x)=\gb$ and $\gb<\gd(y)<(1+a)\gb$, inequality
 \eqref{Har1} holds. Moreover, by Theorem \ref{G-estA},
$$\rec{c}|x-y|^{2-N}\le G_V(x,y)\le c |x-y|^{2-N}.$$
Therefore, as in \eqref{temp2}, we obtain
\begin{equation}\label{temp3'}\BAL
I_{1,3}(\gb)
&\sim \rec{\gb} \int_{\Gs_\gb}\int_{\gb\le \gd(y)}|x-y|^{2-N}\gc_{a\gb}(|x-y|) \Phi_V(y) d\tau(y)dS_x \\
&\sim\rec{\gb}\int_{\gb\le \gd(y)} \int_{\Gs_\gb}|x-y|^{2-N}\gc_{a\gb}(|x-y|)dS_x \Phi_V(y) d\tau(y)\\
&\lesssim \int_{\gb\le \gd(y)}\Phi_V d\tau
\EAL\end{equation}

We also have a (partial) estimate from below.

If $y$ is a point \sth $\gb\leq \gd(y)< \frac{3a}{2}\gb$ then
 $B_{a\gb}(y)\cap \Gs_\gb$ contains an $(N-1)$ dimensional ball of radius $\gb/2$ and consequently there exists a constant $c_3(a)>0$ \sth
 $$\int_{\Gs_\gb}|x-y|^{2-N}\gc_{a\gb}(|x-y|)dS_x>c_3.$$
 Therefore
\begin{equation}\label{est2.3*}
c_3\int_{\gb\leq \gd(y)< \frac{3a}{2}\gb}\Phi_V d\tau \leq I_{1,3}(\gb)
\end{equation}

Inequalities \eqref{est2.3*}, \eqref{est2.3} and \eqref{est2.3'}, imply the lower estimate in \eqref{I1}. The upper estimate follows from \eqref{temp3'}, \eqref{est2.3} and \eqref{est2.3'}.

\qed

\begin{lemma}\label{L-I2}
	Assume that conditions (A1), (A2) and (C1) hold.

		Then there exists $C>0$ \sth  for every $\tau\in \GTM_+(\Gw;\Phi_V)$ and $\gb\in (0,\gb_0)$,
\begin{equation}\label{Gtau2}\BAL
I_{2,\gl}(\gb):=&\rec{\gb^\gl}\int_{\Gs_\gb}\frac{\Phi_V(x)}{\gb}\int_\Gw G_V(x,y) (1-\gc_{a\gb}(|x-y|))d\tau(y) dS_x\\ \leq& C\int_{\Gw}\Phi_V d\tau,
	\EAL	\end{equation}
	where $\gl := 2(\ga^*-\ga)$ when $\ga\neq 1/2$. By \eqref{ga*} $-1<\gl\leq 0$. 
		If $\ga=1/2$ then \eqref{Gtau2} holds with $\gl= 2(\ga^*-\ga-\ge)$ where $\ge>0$ is sufficiently small so that $-1<\gl$. 
\end{lemma}

\proof
Since $\gd(x)=\gb$ and $|x-y|\geq a\gb$,
$$a\inf(\gd(x),\gd(y))\le |x-y|.$$
Therefore, by Theorem \ref{G-estB}, \eqref{est.G3} holds. As before we choose
$x_y= \eta+|x-y|\mbfn_{\eta}$ where $\eta:=\gs(x)$ is the nearest point to $x$ on $\bdw$.  Thus $x$ and $x_y$ are on a normal to $\bdw$ and $|x-y|=\gd(x_y)\geq a\gd(x)$.

By assumption (C1),
$$\frac{\Phi_V(x)}{\Phi_V(x_y)} \leq c(a)\frac{\gb^{\ga^*}}{|x-y|^{\ga}}.$$
Therefore, by \eqref{est.G3},
\begin{equation}\label{est.I2}\BAL
I_{2,0}&=\rec{\gb}\int_{\Gs_\gb}\Phi_V(x)\int_\Gw G_V(x,y) (1-\gc_{a\gb}(|x-y|))d\tau(y) dS_x\\
 &\lesssim \rec{\gb}\int_{\Gs_\gb\cap [|x-y|>a\gb]}\frac{\Phi_V(x)^2}{\Phi_V(x_y)^2}|x-y|^{2-N}\int_{\Gw}\Phi_V(y)d\tau(y)\, dS_x\\
&\lesssim \gb^{2\ga^{*}-1} \int_{\Gs_\gb\cap[|x-y|>a\gb]}|x-y|^{2-N-2\ga} \int_{\Gw}\Phi_V(y)d\tau(y)\, dS_x\\
&\lesssim \gb^{2\ga^*-1} \int_{a\gb}^{1}r^{-2\ga}dr \int_{\Gw}\Phi_V(y)d\tau(y)\\
&\lesssim  \gb^{2\ga^*-2\ga}\int_\Gw \Phi_V(y)d\tau(y).
\EAL\end{equation}
Here we assumed that $\ga\neq 1/2$. In this case \eqref{Gtau2} holds with $\gl= 2\ga^*-2\ga$. If $\ga=1/2$ we replace it by $\ga+\ge$ with $\ge>0$ sufficiently small so that $\gl=2(\ga^*-\ga-\ge)>-1$.
\qed

%
\vskip 4mm

\noindent\textbf{Proof of Theorem \ref{G-I}.}\hskip 2mm By Lemma \ref{L-I1},
$$	\rec{c}\int_{\Gw_{3a\gb/2}}\Phi_V d\tau \leq I_1(\gb)$$
for every $\gb<\gb_0/3a$. Therefore for every $\gb_0/6a<\gb<\gb_0/3a$,
$$	\rec{c}\int_{\Gw_{\gb_0/4}}\Phi_V d\tau \leq I_1(\gb).$$
 More precisely, there exists a constant $c^*$ depending only on $a$,   $\bar a$, $\gb_0$, $\kappa$ \sth
\begin{equation}\label{I1beta}
 c^*\int_{\Gw_{\gb_0/4}}\Phi_V\,d\tau\leq
\int_{[\frac{\gb_0}{6a}<\gd<\frac{\gb_0}{3a}]}\frac{\Phi_V}{\Phi_0}\BBG_V[\tau] dx.
\end{equation}
A suitable constant is given by $c^*= (\inf_{\Gw_{\gb_0/4}}H)^{-1} \frac{\gb_0}{6ac}$ where  $H$ is  the Jacobian  of the transformation from Euclidean  coordinates  to flow coordinates $(\gd,\gs)$. It is known that $H(x)\to 1$ as $\gd(x)\to 0$.

For $a=16(2+\kappa)^2$, \eqref{I1beta} implies \eqref{G-I1}.

Put $\tau'=\tau \chr{[\gd\geq \gb_0/4]}$. Then, using Theorem \ref{G-estA} we obtain,
\begin{equation}\label{I1beta'}\BAL
&\int_{\Gw}\frac{\Phi_V}{\Phi_0} \int_\Gw G_V(x,y)d\tau'(y) dx =
\int_{\Gw} \int_\Gw \frac{\Phi_V}{\Phi_0} G_V(x,y) dx d\tau'(y)\geq\\
&\int_{[\gd(y)\geq \gb_0/4]} \int_{|x-y|<\gb_0/8} \frac{\Phi_V}{\Phi_0} G_V(x,y) dx d\tau(y)\geq\\
c_1&\int_{[\gd(y)\geq \gb_0/4]} \int_{|x-y|<\gb_0/8} \frac{\Phi_V}{\Phi_0} |x-y|^{2-N} dx d\tau(y)\geq\\
c_2&\int_{[\gd(y)\geq \gb_0/4]} \int_{|x-y|<\gb_0/8} |x-y|^{2-N} dx d\tau(y)\geq\\
 c_3 &\int_{[\gd\geq \gb_0/4]} d\tau\geq c_4 \int_{[\gd\geq \gb_0/4]}\Phi_V d\tau,
\EAL\end{equation}
the constants depending only on  $a$,   $\bar a$, $\gb_0$, $\kappa$. 
Combining \eqref{I1beta} and \eqref{I1beta'} we obtain \eqref{G-I1'}.
\qed
\vskip 3mm

\noindent\textbf{Proof of Theorem \ref{G-II}.} \hskip 2mm By \eqref{I1} and \eqref{Gtau2}
$$I_1(\gb)\leq c_1\int_\Gw \Phi_V\,d\tau,\q I_{2,\gl}(\gb)\leq
 c_2\int_\Gw \Phi_V\,d\tau$$
 for every $\gb \in (0, r_1)$ where $r_1:=\gb_0/48(2+\kappa)^2$. Therefore
 $$\int_{\Gs_\gb}\frac{\Phi_V}{\gb} \BBG_V[\tau]dx \leq I_1(\gb) + \gb^{\gl}I_{2,\gl}(\gb)\leq
  c\max(1, \gb^{\gl})\int_\Gw \Phi_V\,d\tau $$
with $\gl$ as in \eqref{Gtau2}. Consequently, integrating over $\gb$ in $(0, r_1)$,
\begin{equation}\label{temp4}
\int_{\Gw_{r_1}} \frac{\Phi_V}{\Phi_0} \BBG_V[\tau]dx \leq C_1 \int_\Gw \Phi_V\,d\tau
\end{equation}
where $C_1$ depends on $\bar a$, $\gb_0, \kappa$, $\ga^*$, $\ga$.

Therefore, to obtain  \eqref{G-I2}, it remains to show that
\begin{equation}\label{temp5}
\int_{D_{r_1}} \frac{\Phi_V}{\Phi_0} \BBG_V[\tau]dx \leq C_2 \int_\Gw \Phi_V\,d\tau
\end{equation}
with $C_2$ depending on the parameters mentioned above.


Let $r_2=r_1/2$ and write,
$$\BAL
\int_{D_{r_1}} \frac{\Phi_V}{\Phi_0} \BBG_V[\tau]dx=
&\int_{D_{r_1}} \frac{\Phi_V}{\Phi_0}\int_{D_{r_2}} G_V(x,y)d\tau(y) dx\\
+&\int_{D_{r_1}} \frac{\Phi_V}{\Phi_0}\int_{\Gw\sms D_{r_2}} G_V(x,y)d\tau(y) dx =:J_1+J_2.
\EAL$$
In $J_2$, $x\in D_{r_1}$ and $y\in \Gw\sms D_{r_2}$. Therefore
$G_V(x,y) \sim \Phi_V(y)$. Consequently
\begin{equation}\label{temp6.1}
J_2 \lesssim \int_{\Gw\sms D_{r_2}}\Phi_V\, d\tau.
\end{equation}
In $J_1$ $x,y\in D_{r_2}$. Therefore, by \cite{Herve}, $G_V(x,y) \sim |x-y|^{2-N}$. Consequently 
\begin{equation}\label{temp6.2}
J_1 \lesssim \int_{D_{r_1}}\int_{D_{r_2}}|x-y|^{2-N} d\tau(y)dx \lesssim
\tau(D_{r_2})\lesssim \int_{D_{r_2}}\Phi_V d\tau.
\end{equation}
Combining these inequalities we obtain \eqref{temp5}.
\qed

\vskip 3mm

\section{A few examples}

We discuss some families of potentials satisfying conditions
(A1), (A2), (C1) required in theorems 3.1 -- 3.3. 

  Denote by $V_F$ a potential of the form 
$$V_F=\rec{\gd_F^2}, \qtxt{$F\sbs\bdw$ a compact set}, \q \gd_F(x)= \dist (x,F).$$

{\bf I.}
Obviously a potential $\gg V_F$ \sth 
$\gg < c_H(V_F)$ satisfies (A1) and (A2).\\ [2mm] 
\textit{If $\gg\in [0,1/4)$ then $\gg V_F$ satisfies condition (C1).} \\[3mm]
 Indeed in this case, 
$$0\leq \gg V_F\leq \gg V_{\bdw}.$$
Consequently, 
$$G_0\leq G_{\gg V_F} \leq G_{\gg V_{\bdw}}$$ 
which implies
\begin{equation}\label{Phi_VF}
c_1\gd\leq  \Phi_{\gg V_F} \leq c_2 \gd^\ga
\end{equation}
for some $\ga\in (\rec{2},1)$. Thus $\gg V_F$ satisfies condition (C1).\\

Similarly, let $V=\gg V_{\bdw}$ where $\gg$ is a bounded measurable function in $\Gw$ \sth 
$$a_1\leq \gg \leq a_2<1/4.$$
Put $\ga_i:=\rec{2}+\sqrt{\rec{4}-a_i}$, $i=1,2$ and assume that
\begin{equation}\label{ga1,2}
 0\leq \ga_1-\ga_2<1/2.
\end{equation}
Then $1/2< \ga_2\leq \ga_1, $ and $$\gd^{\ga_1} \lesssim\Phi_V \lesssim \gd^{\ga_2}.$$ 
Consequently $V$ satisfies condition (C1).\\ [2mm]

{\bf II.} Let $V$ be a potential satisfying (A1), (A2) \sth
\begin{equation}\label{gags}
c_1  \gd^{\ga_1} \gd_F^\gs \leq \Phi_V \leq c_2 \gd^{\ga_2} \gd_F^\gs \quad  *
\end{equation} 
where $\ga_1, \ga_2$ satisfy \eqref{ga1,2}.

Such an estimate holds, for instance, when $V=\gg V_F$, $F\sbs \bdw$ is a smooth k-dimensional manifold without boundary,  $0\leq k\leq N-2$ and $\gg<c_H(V)$.
In this case, it is known that $c_H(V)\leq \rec{4}(N-k)^2$ and, by [MT2], 
\begin{equation}
\Phi_V \sim \gd \gd_F^\gs, \q \gs=\rec{2}\big(k-N+\sqrt{(N-k)^2-4\gg} \,\big).
\end{equation}
Note that $\gs>0$ when $\gg<0$ and $\gs<0$ when $0<\gg<c_H(V)$.\footnote
{Using this example it is possible to construct a family of potentials $V$ for which $\Phi_V$ satisfies \eqref{gags} or related estimates.}

Next we prove,
\begin{lemma}\label{V_C1}
	Let $V=\gg V_F$ where $F\sbs \bdw$ is a compact set and $\gg < c_H(V_F)$. Assume that \eqref{gags} and \eqref{ga1,2} hold. Then $V$ satisfies condition (C1).
\end{lemma}
\Remark Without loss of generality we may assume that $\ga_1+\gs\neq 1/2$. Otherwise we replace $\ga_1$ by $\ga_1+\ge$ with $\ge>0$ \sth $\ga_1+\ge-\ga_1<1/2$.

\proof Let $x,y\in \Gw_{\gb_0}$ and suppose that $\gb_0\geq |x-y|>a\gd(x)$ for some $a\geq 1+a_0$ where $a_0:=16(2+\kappa)^2$. Obviously $\gd(x)\leq \gd_F(x)$; we consider the following cases separately:

$\,$\hskip 2mm (i) $a|x-y|\leq \gd_F(x)$, \q (ii)  $|x-y|\geq a\gd_F(x)$, \q (iii) $|x-y| \sim \gd_F(x)$.\\ [2mm]
Let $z:=(\gd(x)+ |x-y|)\mbfn_{\gs(x)}$. Then $\gd(z)\geq |x-y|>a\gd(x)$ and 
\begin{equation}\label{dF,xz}
\BAL
&|\gd_F(x)-\gd_F(z)|\leq |x-z|=|x-y|,\\    
&a\gd(x)\leq |x-y|\leq \gd(z)
\EAL\end{equation}

\noindent(i) In this case  \eqref{dF,xz} implies
$$\gd_F(z)\leq |x-y|+\gd_F(x)\Longrightarrow\gd_F(z)\lesssim \gd_F(x)$$
and
$$\gd_F(x)-|x-y|\leq \gd_F(z) \Longrightarrow (1-\rec{a})\gd_F(x)\leq \gd_F(z).$$
Thus
\begin{equation}
\gd_F(z)\sim  \gd_F(x).
\end{equation}
Therefore, by \eqref{gags},
\begin{equation}\label{case,i}
\Phi_V(x)/ \Phi_V(z)\lesssim \gd(x)^{\ga_2}/\gd(z)^{\ga_1}.
\end{equation}

\noindent(ii) In this case \eqref{dF,xz} and the definition of $z$ yield
\begin{equation}\label{xz2}\BAL
|x-y|&\leq\gd_F(z)\leq (1+\rec{a})|x-y|, \\ \gd_F(x)&\leq \gd_F(z) + |x-y|
\leq 3|x-y|.
\EAL\end{equation} 
If $\gs<0$ then, by \eqref{dF,xz} and \eqref{xz2},
\begin{equation}\label{case,ii}
\Phi_V(x)/ \Phi_V(z)\lesssim \gd(x)^{\ga_2+\gs}/(\gd(z)^{\ga_1}|x-y|^\gs) \lesssim \gd(x)^{\ga_2+\gs}/\gd(z)^{\ga_1+\gs}.
\end{equation}
If $\gs>0$,
\begin{equation}\label{case,ii'}
\Phi_V(x)/ \Phi_V(z)\lesssim \gd(x)^{\ga_2}|x-y|^\gs/\gd(z)^{\ga_1}|x-y|^\gs.
\end{equation}
Thus, for every $\gs$, \eqref{case,i} holds.
\noindent(iii) In this case,  
$$
|x-y|\leq\gd_F(z)\lesssim |x-y|.
$$ 
Consequently $\gd_F(x) \sim \gd_F(z)$ and \eqref{case,i} follows.
\qed

\section{The $L_V$ boundary trace and positive $L_V$ superharmonic  and subharmonic functions}\label{HM}

\subsection{The boundary trace}
 Assume that $V$ satisfies conditions (A1), (A2).

Let $D\Subset \Gw$ be a Lipschitz domain,  denote by $P_V^D$  the Poisson kernel of $L_V$ in $D$ and by $\gw_V^{x_0,D}$ the harmonic measure of $L_V$
in $D$ relative to a fixed reference point $x_0\in D$.
Then,
\begin{equation}\label{hr-m1}
d \gw_V^{x_0,D} = P_V^D(x_0,\cdot)dS \qtxt{on } \; \prt D.
\end{equation}

Let $\{D_n\}$ be a uniformly Lipschitz exhaustion of $\Gw$. It is well known that if  $u$ is a positive $\Gd$-harmonic function then
\begin{equation}\label{gw0}
u\lfloor_{_{\prt D_n}} d\gw_0^{x_0, D_n} \rightharpoonup \nu
\end{equation}
where $\nu\in \GTM(\bdw)$ is the boundary trace of $u$ and $\rightharpoonup$ indicates weak convergence in measure.  
In \cite[Definition 3.6]{MV_Lip}, \eqref{gw0} was used as a definition of boundary trace for  solutions of certain semilinear equations with absorption. 
In this spirit we define,
\begin{definition}\label{ntr}  	A non-negative Borel function $u$ defined in $\Gw$ has an $L_V$ \emph{boundary trace} $\nu\in \GTM(\bdw)$ if
	\begin{equation}\label{ntr1}
	\lim_{n\to\infty}\int_{\prt D_n}hud\gw_V^{x_0,D_n}= \int_{\bdw} hd\nu \forevery h\in C(\bar\Gw),
	\end{equation}
for every uniformly Lipschitz exhaustion $\{D_n\}$ of $\Gw$.
	The $L^V$ trace will be denoted by $\tr_V(u)$.
		Here  we assume that $D_{\gb_0}\sbs D_1$ and $x_0 \in D_{\gb_0}$.
	\end{definition}

\begin{lemma}\label{weak_con} If $\{D_n\}$ is a uniformly Lipschitz exhaustion of $\Gw$ then, for every positive $L_V$ harmonic function $u=\BBK_V[\nu]$,
	\begin{equation}\label{hr-meas}
	\lim_{n\to \infty}\int_{\prt D_n}hu\, d\gw_V^{x_0,D_n} = \int_{\bdw}h\,d\nu
	\forevery h\in C(\bar \Gw).
	\end{equation}
\end{lemma}
\noindent\textit{Remark.}\hskip 2mm The proof is similar to that of \cite[Lemma 2.2]{MV_Lip}. We omit details.

\begin{lemma}\label{trKV}
	Assume (A1) and (A2). Then
	\begin{equation}\label{trKV1}\BAL
(a)\q	&\tr_V(\BBK_V[\nu])=\nu &\forevery \nu\in \GTM_+(\bdw)\\
(b)\q 	&\tr_V(\BBG_V[\tau])=0  &\forevery \tau\in \GTM_+(\Gw;\Phi_V).
	\EAL\end{equation}
\end{lemma}

\proof
(a) is a restatement of Lemma \ref{weak_con}. 

(b) follows from the fact that $\BBG_V[\tau]$ is an $L_V$ potential , i.e., it does not dominate any positive $L_V$ harmonic function (see \cite{An-SLN}).

 Let $\{D_n\}$ be a Lipschitz exaustion of $\Gw$ and let
   $$u_n(x) = \int_{\prt D_n} \BBG_V[\tau] P_V^{D_n}(x,\gx)dS_\gx \forevery x\in D_N.$$
 Then $u_n$ is $L_V$ harmonic in $D_n$ and it is dominated by $G_V[\tau]$ which is \LVsup. Hence 
 $u=\lim u_n$ is an $L_V$ harmonic function dominated by $\BBG_V[\tau]$. By the previous observation, $u\equiv 0$. 
 
 \qed

%
%

%
%
%

\subsection{Estimates of positive \LVsup and \LVsub functions.}
\begin{proposition}\label{sub-sup} Assume (A1), (A2).\\
	(i) If  $u$ is a positive \LVsup function then 
	\begin{equation}\label{sup1}\BAL
	-L_Vu=\tau \in &\GTM_+(\Gw;\Phi_V) 
	\EAL\end{equation}
	and there exists a non-negative measure $\nu\in \GTM(\bdw)$ \sth
	\begin{equation}\label{sup-rep}
	u=\BBG_V[\tau]+\BBK_V[\nu].
	\end{equation}
	\noindent(ii) 	Let $u$ be a non-negative \LVsub function  and  $\tau:=L_V u$. (Thus $\tau\geq 0$ is a Radon measure.) Then, 
	\begin{equation}\label{tau_good}
\tau  \in \GTM(\Gw;\Phi_V) \Longleftrightarrow \text{$u$ has an $L_V$ boundary trace}
	\end{equation}
If 	$\nu:=\tr_V u$ exists then,
%
%
%
%
	\begin{equation}\label{sub-rep}
	u+\BBG_V[\tau]=\BBK_V[\nu].
	\end{equation}

\end{proposition}

\proof
(i) This statement is an immediate consequence of the Riesz decomposition lemma and the fact that $u$ is an  $L_V$ potential if and only if $u=\BBG_V[\gs]$ for some $\gs \in \GTM_+(\Gw;\Phi_V)$.

\noindent(ii) If $\tau \in \GTM(\Gw;\Phi_V)$ then $u+\BBG_V[\tau]$ is positive $L_V$ harmonic. By the Representation Theorem  $\exists\; \nu\in \GTM(\bdw)$ \sth $u+\BBG_V[\tau]=\BBK_V[\nu]$. 
By Lemma \ref{trKV} $\nu=\tr_Vu$.

Conversely suppose that $\nu=\tr_Vu$ esists. 


Let $\{D_n\}$ be a Lipschitz exaustion of $\Gw$. The function $u_n=u\chr{D_n}$ satisfies,
\begin{equation}\label{bvpn}
u_n+ \BBG_V^{D_n}[\tau\chr{D_n}]= \int_{\prt D_n} u P_V^{D_n}(x,\gx)dS_\gx.
\end{equation}

Since $\nu:= \tr_V u$ exists, the last term converges at the point $x_0$ to $\nu$. Consequently $\{ \BBG_V^{D_n}[\tau\chr{D_n}](x_0)\}$ is bounded and (by the monotone convergence theorem) converges to 
$\BBG_V[\tau](x_0)$. Therefore $\tau \in \GTM(\Gw;\Phi_V)$.

\qed

Combining Proposition \ref{sub-sup} with Theorems \ref{Knu}, \ref{G-I} and \ref{G-II} we obtain the folowing two sided estimates.

\begin{theorem}\label{LVsup} Assume (A1), (A2), (C1).
	
	(i)	Let  $u$ be a positive \LVsup function and let $\tau$, $\nu$ be as in Proposition \ref{sub-sup}. Then there exists a constant $C$ depending only on $\bar a$, $\ga^*$, $\ga$, $\gb_0$, $\Gw$  \sth
	\begin{equation}\label{LVsup1}
	\rec{C}(\int_{\Gw} \Phi_V\,d\tau + \norm{\nu}) \leq \int_\Gw \frac{\Phi_V}{\Phi_0}\,udx\leq C(\int_\Gw \Phi_V\,d\tau + \norm{\nu}).
	\end{equation}

	(ii)	Let  $u$ be a positive \LVsub function and 
	assume that  
	\begin{equation}\label{good_sub}
	\tau:=L_V u  \in \GTM(\Gw;\Phi_V).
	\end{equation}
	Then \eqref{sub-rep} holds and there exists a constant $C$ depending only on $\bar a$, $\ga^*$, $\ga$, $\gb_0$, $\Gw$ \sth
	
	\begin{equation}\label{LVsub1}
	C^{-2} \norm{\nu} \leq  \int_\Gw \frac{\Phi_V}{\Phi_0}\,u\,dx + \int_{\Gw} \Phi_V d\tau \leq C^2 \norm{\nu}
	\end{equation}
\end{theorem}
\proof Inequality \eqref{LVsup1} is an immediate consequence  of \eqref{sup-rep} and the estimates stated in Theorems \ref{Knu}, \ref{G-I} and \ref{G-II}.

For part (ii) we use \eqref{sub-rep} and the above mentioned estimates to obtain:
\begin{equation}\label{LVsub2}\BAL
\int_\Gw \frac{\Phi_V}{\Phi_0}\,udx + \rec{C}\int_{\Gw} \Phi_V d\tau \leq 	\int_\Gw \frac{\Phi_V}{\Phi_0}\,(u +\BBG[\tau])dx\leq C \norm{\nu}, \\
\rec{C}  \norm{\nu}\leq \int_\Gw \frac{\Phi_V}{\Phi_0}\,(u +\BBG[\tau])dx \leq  \int_\Gw \frac{\Phi_V}{\Phi_0}\,u dx + C\int_{\Gw} \Phi_V d\tau.		
\EAL	\end{equation}	
\qed

\vskip 20mm

\small{Address: Department of Mathematics, Technion, Haifa, Israel

E-mail: marcusm@math.technion.ac.il}

\newpage

\end{document}